\documentclass[a4paper,12pt]{article}

\usepackage{amssymb}
\usepackage{amsmath}
\usepackage{amsfonts}                                                           
\usepackage{mathrsfs} 
\usepackage{amsthm}

\usepackage[table]{xcolor}  

\usepackage{multirow}
\usepackage{colortbl}

\definecolor{TableHeadGray}{gray}{.8}

\DeclareFontFamily{OT1}{pzc}{}
\DeclareFontShape{OT1}{pzc}{m}{it}{<-> s * [1.200] pzcmi7t}{}
\DeclareMathAlphabet{\mathpzc}{OT1}{pzc}{m}{it} 

\newcommand{\Icg}[2]{\mathrm{ICG}({#2},{#1})}   
\newcommand{\Ene}[2]{\mathpzc{E}({#2},{#1})}   
  
\newcommand{\Emin}[1]{\mathpzc{E}_{\mathrm{min}}({#1})}   
\newcommand{\Emax}[1]{\mathpzc{E}_{\mathrm{max}}({#1})}

\def \C{\hbox{\sf \rlap{\kern.25em \vrule width.1em height1.6ex depth-.1ex}C}}
\def \D{{\sf I\kern-1.5ptD\,}}
\def \K{{\sf I\kern-1.5ptK\,}}
\def \N{{\sf I\kern-1.5ptN\,}}
\def \P{{\sf I\kern-1.5ptP\,}}
\def \Q{\hbox{\sf \rlap{\kern.25em \vrule width.1em height1.6ex depth-.1ex}Q}}
\def \R{{\sf I\kern-1.5ptR\,}}
\def \T{{\sf T\kern-6.5ptT\,}}
\def \Z{{\sf Z\kern-5.0ptZ\,}}

\def\begitem#1 {\bigskip\pagebreak[1]%
     \refstepcounter{subsection}{\nopagebreak[4]%
     \thesubsection\hskip 0.5truecm}
     {\sc#1}\hskip 1pt.\nopagebreak[4]\par\nopagebreak[4]%
      \begin{enumerate}\rm\nopagebreak[4]}
\def\BEGITEMK#1 #2{\bigskip\pagebreak[1]%
      \refstepcounter{subsection}{\nopagebreak[4]%
     \thesubsection\hskip 0.5truecm}\nopagebreak[4]
     {\bf#1}\hskip 1pt.\nopagebreak[4]\par\nopagebreak[4]%
     \medskip\nopagebreak[4]\rm#2\nopagebreak[4]%
     \begin{enumerate}\nopagebreak[4]\rm}
\def\enditem{\end{enumerate}}

\renewcommand{\qed}{\hfill \ensuremath{\Box}\bigskip}

\newtheorem{theorem}{Theorem}[section]

\newtheorem{lemma}{Lemma}[section]
\newtheorem{proposition}{Proposition}[section]
\theoremstyle{definition}
\newtheorem*{definition}{Definition}

\newtheorem{example}{Example}[section]

\begin{document}


\title{The exact maximal energy of integral circulant graphs with prime power order}

\author{{\bf J.W.~Sander\footnote{Corresponding author}}\\
Institut f\"ur Mathematik und Angewandte Informatik,
Universit\"at Hildesheim,\\
D-31141 Hildesheim, Germany\\
{\tt sander@imai.uni-hildesheim.de}\\[.1in]
and\\[.1in]
{\bf T.~Sander}\\ 
Fakult\"at f\"ur Informatik,
Ostfalia Hochschule f\"ur angewandte Wissenschaften,\\
D-38302 Wolfenb\"uttel, Germany\\
{\tt t.sander@ostfalia.de}}

\maketitle

\begin{abstract}

The energy of a graph was introduced by {\sc Gutman} in 1978 as the sum of the absolute values of the eigenvalues of its adjacency matrix. We study the energy of integral circulant graphs, also called gcd graphs,
which can be characterized by their vertex count $n$ and a set $\cal D$ of divisors of $n$ in such a way that they have vertex set $\mathbb{Z}/n\mathbb{Z}$ and edge set $\{\{a,b\}:\, a,b\in\mathbb{Z}/n\mathbb{Z},\, \gcd(a-b,n)\in {\cal D}\}$. 

Given an arbitrary prime power $p^s$, we determine all divisor sets maximising the energy of an integral circulant graph of order $p^s$. This enables us to compute the maximal energy $\Emax{p^s}$ among all integral circulant graphs of order $p^s$. 
\end{abstract}

{\bf 2010 Mathematics Subject Classification:} Primary 05C50, Secondary 15A18
\\

{\bf Keywords:} Cayley graphs, integral graphs, circulant graphs, gcd graphs, graph energy
\\

%
%

%
%

\section{Introduction}\label{intro}
%

Given a positive integer $n$ and a set ${\cal D}$ of positive divisors of $n$, the \textit{integral circulant graph} $\Icg{\cal D}{n}$ is  
defined as the graph having the residue class ring  $\mathbb{Z}/n\mathbb{Z}$
as vertex set and  $\{ \{a,b\}:~a,b\in \mathbb{Z}/n\mathbb{Z},~ \gcd(a-b,n)\in {\cal D}\}$ as edge set. These graphs are also known as \textit{gcd graphs} in the literature.
For $\vert {\cal D}\vert =1$ we obtain the subclass of so-called \textit{unitary Cayley graphs}.
We consider only loopfree gcd graphs, i.e. we require $n\notin {\cal D}$. Moreover, we note that
$\Icg{\cal D}{n}$ with ${\cal D} = \{d_1,\ldots,d_r\}$ is connected if and only if  $\gcd(n,d_1,\ldots,d_r)=1$ (cf.~\cite{SO}). For a prime power $n=p^s$ this is equivalent to $1\in{\cal D}$.

In recent years, quite a few structural properties of integral circulant graphs have been brought to light
(cf. \cite{DEJ}, \cite{BER}, \cite{SO}, \cite{KLO}, \cite{AKH}, 
\cite{BAS}, \cite{KLO2}, \cite{DRO}, \cite{BEA}, \cite{BAS2}).
Some emphasis has lately been placed on researching the energy of integral circulant graphs 
(see \cite{SHP}, \cite{ILI}, \cite{RAM}, \cite{BAS3}, \cite{PET}, \cite{SA1}, \cite{SA2}, \cite{KIA}).
The \textit{energy} $E(G)$ of a graph $G$ on $n$ vertices is defined as
\[
     E(G)=\sum_{\lambda \in {\rm Spec}(G)} \vert \lambda \vert \quad =\; \sum_{i=1}^n \vert\lambda_i\vert,
\]
where the \textit{spectrum} ${\rm Spec}(G) =\{ \lambda_1,\ldots,\lambda_n\}$ of $G$ denotes the set of 
eigenvalues $\lambda_i$ of the adjacency matrix of $G$, counted with multiplicities. Observe that an undirected graph has real spectrum, since its adjacency matrix is symmetric.

Let us abbreviate $\Ene{\cal D}{n}=E(\Icg{\cal D}{n})$. Since $\Icg{\cal D}{n})$ is an integral graph, we have ${\rm Spec}(\Icg{\cal D}{n}) \subset \mathbb{Z}$, hence $\Ene{\cal D}{n}$ is an integer. 
Given a positive integer $n$, it is most desirable to determine
\begin{align*}
\Emin{n} &:= \min\,\{ \Ene{\cal D}{n}:\;\; {\cal D}\subseteq \{1\le d<n:\; d\mid n\} \} \\[-.2in]
\intertext{and} \\[-.5in]
\Emax{n} &:= \max\,\{ \Ene{\cal D}{n}:\;\; {\cal D}\subseteq \{1\le d<n:\; d\mid n\} \},
\end{align*}
being of particular interest with respect to the question of hyperenergeticity or hypoenergeticity of certain classes of graphs.

Consider a prime power $n=p^s$ and an arbitrary divisor set ${\cal D} = \{ p^{a_1}, p^{a_2} , \ldots , p^{a_r}\}$ with exponents $0\le a_1 < \ldots < a_r \le s-1$. 
In \cite{SA1}, Theorem 2.1, the authors proved that
\begin{equation}
\Ene{{\cal D}}{p^s} = 2(p-1)p^{s-1}\left(r -(p-1)h_{p,r}(a_1,\ldots,a_r) \right),  \label{ft3}
\end{equation}
where
\begin{equation}
h_{p,r}(x) = h_{p,r}(x_1,\ldots,x_r) := \sum_{k=1}^{r-1}\sum_{i=k+1}^r \,\frac{1}{p^{x_i-x_k}}   \label{basishp}
\end{equation}
for $x=(x_1, \ldots, x_r) \in \mathbb{R}^r$. Observe that $h_{p,r}$ has the symmetry property
\begin{equation}
 h_{p,r}(s-1-a_r,\ldots,s-1-a_1) = h_{p,r}(a_1,\ldots,a_r) \label{symmhp}
 \end{equation} 
for all integral exponents $0\leq a_1 < a_2 < \ldots < a_{r-1} < a_r \leq s-1$.

A rather straightforward consequence of (\ref{ft3}) is that 
$$ \Emin{p^s} = 2(p-1)p^{s-1},$$
and this minimal energy is attained precisely for the singleton divisor sets ${\cal D}=\{p^t\}$ with $0\le t \le s-1$ (cf. \cite{SA1}, Theorem 3.1).

Divisor sets $\cal D$ producing graphs with maximal energy 
$ \Emax{p^s}$ were studied for the first time in \cite{SA2}. 
Since, by the aforementioned result, singleton divisor sets generate integral circulant graphs having minimal energy $\Emin{p^s}$, minimisers of some $h_{p,r}$ producing integral circulant graphs with maximal energy 
$\Emax{p^s}$ necessarily satisfy $r\geq 2$ if $s\geq 2$; note that for $s=1$ there is only one possible divisor set, namely ${\cal D}=\{1\}$. For that reason, we may henceforth assume w.l.o.g. that $2\leq r\leq s$. 
Furthermore, any minimiser $0\leq a_1<a_2<\ldots<a_r\leq s-1$ must have the entries $a_1=0$ and $a_r=s-1$, because otherwise replacing $a_1$ by a smaller integer or $a_r$ by a larger one, respectively, would obviously summandwise lessen the value of $h_{p,r}$ as defined in (\ref{basishp}).  Accordingly, any minimiser $a=(a_1,\ldots,a_r)$ of $h_{p,r}$ lies in the set
\[ A(s,r) := \{(a_1,\ldots,a_r)\in \mathbb{Z}^r:\; 0=a_1<a_2<\ldots <a_{r-1}<a_r=s-1\},\]
and such an $a$ is called an \emph{admissible} exponent tuple. Hence $1=p^{a_1}\in {\cal D}$ for all $a\in A(s,r)$, and by one of our introductory remarks we incidentally have that $\Emax{p^s}$ can only be generated by connected gcd graphs.

Our strategy in \cite{SA2} for finding these maximising divisor sets was to start by fixing $r$ and discover exponent tuples $(a_1,\ldots,a_r)$ minimising $h_{p,r}$, at least approximately. 
Applying methods from convex optimisation it was shown that, for fixed $s$ and $r$, the function $h_{p,r}$ becomes minimal if 
$0=a_1<a_2<\ldots<a_{r-1}<a_r=s-1$
are chosen in nearly equidistant position (\cite{SA2}, Corollary 3.2).
Inserting the corresponding approximations $h_{p,r}(a_1,\ldots,a_r)$ into (\ref{ft3}) and varying $r$,  the order of magnitude of $\Emax{p^s}$ was determined in the sense that      
explicit upper and lower bounds for $\Emax{p^s}$ were given, which differed roughly by a factor $2$. More precisely, it was shown in \cite{SA2}, Theorem 4.2 that 
\[
   \left(1-\frac{\log\log p}{\log p}\right)(p-1)p^{s-1}(s-1)\leq  \Emax{p^s}  \leq 2\left(1-\frac{\log\log p}{\log p}\right)(p-1)p^{s-1}s
\]
for all primes $p\geq 17$. Numerical examples suggested that $\Emax{p^s}$ approached the lower bound with increasing $p$, i.e. we had conjecturally $\Emax{p^s} \approx sp^s$.

In \cite{SA3} the authors used combinatorial instead of analytic arguments to refine the earlier approximative results on minimisers of $h_{p,r}$ for fixed values of $r$. This shed more light on the structure of these minimisers, which can be regarded as the result of a repeated balancing process. In several cases this process allowed us to obtain accurate results after only one or two balancing steps.

In this paper, we completely settle the problem to find all divisor sets maximising the energy of an integral circulant graph of prime power order, i.e. for any prime power $p^s$ we shall explicitly determine all divisor sets $\cal D$ satisfying $\Ene{\cal D}{p^s} = \Emax{p^s}$. At the same time this enables us to compute $\Emax{p^s}$ precisely, thus confirming the conjecture mentioned above. 

For an admissible exponent tuple $a=(a_1,\ldots,a_r)$, we denote by
\[ {\cal D}(a)  := \{p^{a_1}, \ldots, p^{a_r}\}\]
the corresponding divisor set. Our main result is the following

\begin{theorem}{\label{thm}}
Let $p$ be a prime, and let $s$ be a positive integer. 
\begin{itemize}
\item[(i)] If $s$ is odd, then
\begin{equation}  \label{emax1}
\Emax{p^s} = \frac{1}{(p+1)^2}\left((s+1)(p^2-1)p^s + 2(p^{s+1}-1)\right),
\end{equation}
and the only exponent tuple $a$ satisfying $\Ene{{\cal D}(a)}{p^s} = \Emax{p^s}$
is $a=(0,2,4,\ldots,s-3,s-1)\in A(s,r)$ in case $p\geq 3$, while we have the additional maximising tuple $a=(0,1,3,5,\ldots,s-4 ,s-2,s-1)$ in case $p=2$.
\item[(ii)] If $s$ is even, then
\begin{equation}  \label{emax0}
\Emax{p^s} = \frac{1}{(p+1)^2}\left(s(p^2-1)p^s + 2(2p^{s+1}-p^{s-1}+ p^2-p-1)\right),
\end{equation}
and the only exponent tuples $a$ satisfying $\Ene{{\cal D}(a)}{p^s} = \Emax{p^s}$ are $a=(0,2,4,\ldots,s-2,s-1)$ and
$a=(0,1,3,5,\ldots,s-3,s-1)$. 
\end{itemize}
\end{theorem}
\bigskip
We like to add the following remarks with regard to Theorem \ref{thm}:
\begin{itemize}
\item[(a)] 
For computational reasons we prefer to present formulae 
(\ref{emax1}) and (\ref{emax0}) in a compact form. However, being the sum of the moduli of integral eigenvalues, $\Emax{p^s}$ is certainly a positive integer.
In fact, it is an immediate consequence of identity (\ref{ft3}) that $\Emax{p^s}$ is always divisible by $2(p-1)$. This
is in line with the work of {\sc Bapat} and {\sc Pati} \cite{BAP} who showed that
the energy of any graph is never an odd integer (see also {\sc Pirzada} and {\sc Gutman} \cite{PIR}). The divisibility property of $\Emax{p^s}$ becomes obvious by not using geometric sum formulae in (\ref{hpodd}) and (\ref{hpeven}), which leads to
\begin{align*}
\Emax{p^{2m+1}} &=  2(p-1)\left((m+1)p^{2m} -(p-1)\sum_{j=0}^{m-1}(j+1)p^{2j}\right)\\[-.2in]
\intertext{and} \\[-.5in]
\Emax{p^{2m}} &= 2(p-1)\left(m p^{2m-1} -(p-1)\sum_{j=0}^{m-3}(j+1)p^{2j+3} +1 \right)
\end{align*}
instead of (\ref{emax1}) and (\ref{emax0}).
\item[(b)] 
Note that the two different exponent tuples in part (ii) just occur due to the symmetry property (\ref{symmhp}) of the function $h_{p,r}$.
\item[(c)]  A graph $G$ on $n$ vertices is called \textit{hyperenergetic} if its energy is greater than the energy of the complete graph $K_n$, i.e.~if $E(G) > E(K_n) = 2(n-1)$. 
There exist several bounds for the energies of different graph classes. For an arbitrary graph $G$ with $n$ vertices, {\sc Koolen} and {\sc Moulton} \cite{KOO} showed that
\begin{equation} \label{KooMou}
 E(G) \leq \frac{n}{2} \left( \sqrt{n} + 1\right).
\end{equation}
{\sc Shparlinksi} \cite{SHP} constructed an infinite family of circulant graphs that asymptotically achieves the upper bound in (\ref{KooMou}). We observe that integral circulant graphs of prime power order $p^s$ with maximal energy, that is $\Emax{p^s}\approx sp^s$ by Theorem \ref{thm}, are hyperenergetic, but do not come close to the bound in (\ref{KooMou}). The graphs $\Icg{\cal D}{p^s}$ with minimal energy, studied in \cite{SA1}, are \textit{hypoenergetic}, which means $\Emin{p^s} < 2(p^s-1)$. The reader finds more on hyperenergeticity as well as hypoenergeticity of gcd graphs in \cite{SA1} and \cite{SA2}.
\end{itemize}

The proof of Theorem \ref{thm} uses an accordion-like compression and expansion procedure, which generates for any given admissible $a^{(0)}\in A(s):= \bigcup_{r=2}^s A(s,r)$ a finite sequence $(a^{(\ell)})_{0\leq \ell \leq m}$, say,  of admissible $a^{(\ell)} \in A(s)$  satisfying 
\[
\Ene{{\cal D}(a^{(0)})}{p^s} < \Ene{{\cal D}(a^{(1)})}{p^s} < \ldots < \Ene{{\cal D}(a^{(m)})}{p^s},
\]
where $a^{(m)}$ is one of the maximising exponent tuples to be found in Theorem \ref{thm}.
Let us point out that, contrary to earlier strategies applied in \cite{SA2} and \cite{SA3}, we now compare the energies related to admissible exponent tuples of 
\;d\;i\;f\;f\;e\;r\;e\;n\;t\,
lengths. We recommend to take a look at the illustrative Example \ref{exseq}. 

With respect to possible generalisations, let us point out some existing obstacles.
So far, no formula comparable with (\ref{ft3}) in terms of simplicity is available for $\Ene{{\cal D}}{n}$ if $n$ is not a prime power. In particular, the energy reveals practically no signs of
multiplicativity in terms of divisors of $n$.
As a singular result it has been observed in
\cite{SA1} that $\Emax{pq}=\Emax{p}\Emax{q}$ for distinct odd primes $p$ and $q$. 
However, in this case it is a straightforward exercise to determine and compare the energies for the four possible divisor sets $\{1\}, \{1,p\}, \{1,q\}, \{1,p,q\}$ by evaluating the formula
\begin{equation*}
\Ene{{\cal D}}{n} = \sum_{k=1}^n \; \left\vert \;\sum_{d\in {\cal D}}\; \mu\left(\frac{n}{(n,kd)}\right) \cdot
\frac{\varphi\left(\frac{n}{d}\right)}{\varphi\left(\frac{n}{(n,kd)}\right)} \;\right\vert,    
\end{equation*}
(cf. \cite{KLO}, Theorem 16) with M\"obius' function $\mu$ and Euler's totient function $\varphi$. The multiplicativity already vanishes for $n=p^2q$ or a product of three distinct primes.
Moreover, no conjectures regarding the structure of the energy maximising divisor sets exist so far,
not even for square-free $n$.
For example, $\Emax{3\cdot 5\cdot 7}=520$ with unique maximiser $\{1, 15, 21, 35\}$ and
$\Emax{2\cdot 3\cdot 5\cdot 7}=1414$ with unique maximiser $\{1, 2, 3, 30, 35, 42, 70, 105\}$, 
while  $\Emax{p}=2(p-1)$ with unique maximiser $\{1\}$ for each prime $p$.

%
%

\section{Comparison of certain admissible exponent tuples} \label{Compare}

Given $2\leq r \leq s$, we define for each $a\in A(s,r)$ its \emph{delta vector}  
\[
\delta_r(a):= (d_1, d_2, \ldots, d_{r-1}) \in \mathbb{N}^{r-1}
\]
by setting $d_j:=a_{j+1}-a_j$ ($1\le j \leq r-1$). Obviously, we have $\sum_{j=1}^{r-1} d_j = s-1$. Thus, introducing
\[
D(s,r):= \{(d_1,\ldots,d_{r-1})\in \mathbb{N}^{r-1}: \; \sum_{j=1}^{r-1} d_j = s-1 \}, 
\] 
the function
\[ 
\delta_r: \left\{ \begin{array}{ccl}
                   A(s,r) &  \longrightarrow & D(s,r) \\
                   (a_1,a_2,\ldots,a_r) & \mapsto & (a_2-a_1,a_3-a_2,\ldots,a_r-a_{r-1}) 
                   \end{array}  
                   \right.
\] 
is 
1--1 
with its inverse 
\[ \quad\quad \quad\quad\;\;\,
\delta_r^{-1}: \left\{ \begin{array}{ccl}
                   D(s,r) &  \longrightarrow & A(s,r) \\
                   (d_1,d_2,\ldots,d_{r-1}) & \mapsto & (0,d_1,d_1+d_2,\ldots, d_1+d_2+\ldots+d_{r-2},s-1). 
                   \end{array}  
                   \right.                  
\]

As an immediate consequence of \cite{SA1}, Theorem 2.1, we have the following observation, which will be used several times in the sequel.
\begin{lemma}  \label{hilfssatz}
Let $p$ be a prime, and let $s\geq 2$ be a fixed integer. Let $d\in D(s,r)$ and $d'\in D(s,r')$ for some
$2\leq r,r' \leq s$, and define $a:=\delta_r^{-1}(d)\in A(s,r)$ and $a':= \delta_{r'}^{-1}(d')\in A(s,r')$. If 
\begin{equation}  \label{hisavorauss}
h_{p,r}(a) - h_{p,r'}(a') > \frac{r-r'}{p-1},
\end{equation}
then  $\Ene{{\cal D}(a')}{p^s} > \Ene{{\cal D}(a)) )}{p^s}$.
\end{lemma}

{\sc Proof. }
It follows from (\ref{ft3}) that
\[
\Ene{{\cal D}(a')}{p^s} - \Ene{{\cal D}(a)}{p^s} = 2(p-1)p^{s-1}\left(r'-r -(p-1)(h_{p,r'}(a') - h_{p,r}(a)\right),
\]
which is positive by (\ref{hisavorauss}).

\qed

We denote by 
$\|d\|_{\infty}:=\max\{d_j:\, 1\leq j \leq r-1\}$ 
the standard maximum norm of $d=(d_1,\ldots, d_{r-1}) \in D(s,r)$.

\begin{proposition}{\label{maxd=3}}
Let $p$ be a prime, and let $s\geq 2$ be a fixed integer. Assume that $d=(d_1,\ldots,d_{r-1})\in D(s,r)$ for some $2\leq r \leq s$ satisfies one of the following two conditions
\begin{itemize}
\item[(i)] $\|d\|_{\infty} \geq 4$,
\item[(ii)] $\|d\|_{\infty} = 3$ and $d_j\geq 2$ for $1\leq j \leq r-1$.
\end{itemize}
Let $1\leq u \leq r-1$ be such that $d_u=\|d\|_{\infty}$, and define 
\[
d' =(d_1',\ldots,d_{r}'):= (d_1,\ldots,d_{u-1},2,d_u-2,d_{u+1},\ldots,d_{r-1}). 
\]
Then $d'\in D(s,r+1)$ and 
$\Ene{{\cal D}(\delta_{r+1}^{-1}(d'))}{p^s} > \Ene{{\cal D}(\delta_r^{-1}(d))}{p^s}$.
\end{proposition}

{\sc Proof. }
By the fact that $\sum d_j =s-1$ for $d\in D(s,r)$, we have 
\[
1\leq u \leq r-1\leq s-\|d\|_{\infty}\leq s-3,
\] 
thus $r+1\leq s-1$ and $\sum d_j' =s-1$, which implies $d'\in D(s,r+1)$. 
For $a=(a_1,\ldots,a_r) := \delta_r^{-1}(d) \in A(s,r)$, we obtain that
$a' = (a_1',\ldots, a_r',a_{r+1}') := \delta_{r+1}^{-1}(d') \in A(s,r+1)$ satisfies
\[
a_j' = \left\{ \begin{array}{ll}
                   a_j &  \mbox{ for $1\leq j \leq u$,} \\
                   a_u+2 &  \mbox{ for $j = u+1$,} \\
                   a_{j-1} &  \mbox{ for $u+2\leq j \leq r+1$.} 
                   \end{array}  
                   \right.                  
\]
By definition of the functions $h_{p,r}$ and $h_{p,r+1}$ in (\ref{basishp}), we have
\[ h_{p,r+1}(a') - h_{p,r}(a) =  \sum_{k=1}^{r}\sum_{i=k+1}^{r+1} \,\frac{1}{p^{a'_i-a'_k}}  -   \sum_{k=1}^{r-1}\sum_{i=k+1}^r \,\frac{1}{p^{a_i-a_k}}.
\]
Since all entries of $a$ are also entries of $a'$, each difference $a_i-a_k$ occurring in the second double sum also occurs in the first one. Hence the corresponding summands cancel out. However, in comparison with $a$ the tuple $a'$ has the additional entry $a'_{u+1}=a_u+2$. Therefore, $h_{p,r+1}(a') - h_{p,r}(a)$ consists of all summands in which $a'_{u+1}$ shows up, thus
\begin{align} \label{diffhp1}
\begin{split}
h_{p,r+1}(a') - h_{p,r}(a) &= 
\sum_{k=1}^u \,\frac{1}{p^{a'_{u+1}-a'_k}}  +   \sum_{i=u+2}^{r+1} \,\frac{1}{p^{a'_i-a'_{u+1}}}   \\
&= \sum_{k=1}^u \,\frac{1}{p^{(a_u +2)-a_k}}  +   \sum_{i=u+2}^{r+1} \,\frac{1}{p^{a_{i-1}-(a_u+2)}}  \\
&= \sum_{k=1}^u \,\frac{1}{p^{a_u-a_k+2}}  +   \sum_{i=u+1}^r \,\frac{1}{p^{a_i-a_u-2}}. 
\end{split}
\end{align}

\underline{Case (i)}: $d_u \geq 4$.\newline
Since $a_1<a_2<\ldots<a_u$, we have $a_u-a_k\geq u-k$ for $k=1,2,\ldots,u$, and since $a_{u+1}-a_u = d_u \geq 4$, we have $a_i-a_u \geq i-u+3$ for $i=u+1,\ldots,r$. By (\ref{diffhp1}), this implies
\begin{align*}  
h_{p,r+1}(a') - h_{p,r}(a) &\leq
\sum_{k=1}^u \,\frac{1}{p^{u-k+2}}  +   \sum_{i=u+1}^r \,\frac{1}{p^{i-u+1}} 
< 2 \sum_{j=2}^{\infty}\frac{1}{p^j} = \frac{2}{p(p-1)} \leq \frac{1}{p-1}
\end{align*}
for each prime $p$. 

\underline{Case (ii)}: $d_u = 3$ and $d_j\geq 2$ for $j\neq u$.\newline
Since $a_{u+1}-a_u = d_u =3$ and $a_{j+1}-a_j=d_j\geq 2$ for all $j$, we have $a_u-a_k\geq 2(u-k)$ for $k=1,2,\ldots,u$ and $a_i-a_u \geq 2(i-u)+1$ for $i=u+1,\ldots,r$. By (\ref{diffhp1}), this implies
\begin{align*}  
\begin{split}
h_{p,r+1}(a') - h_{p,r}(a) &\leq
\sum_{k=1}^u \,\frac{1}{p^{2(u-k+1)}}  +   \sum_{i=u+1}^r \,\frac{1}{p^{2(i-u)-1}} < \sum_{j=1}^{\infty}\frac{1}{p^{2j}} +\frac{1}{p}\sum_{j=0}^{\infty}\frac{1}{p^{2j}}  = \frac{1}{p-1}\,.
\end{split}
\end{align*}
In both cases we have $h_{p,r}(a) - h_{p,r+1}(a') > -\frac{1}{p-1}$, hence Lemma \ref{hilfssatz}
proves our Proposition.

\qed

Before we proceed we introduce a helpful tool to visualise the calculation of differences $h_{p,r}(a)-h_{p,r'}(a')$. Given $a=(a_1,\ldots,a_r)\in A(s,r)$ for some $2\leq r \leq s$, we define the \emph{delta tableau} $\mathfrak{A}=(a_{k,i})_{1\leq k<i\leq r}$ of $a$ as a triangular 
array of integers, corresponding to a strictly upper triangular matrix 
$(a_{k,i})_{1\leq k,i\leq r}$ with $a_{k,i}:= a_i-a_k$ for $1\leq k<i\leq r$ and no entries attributed to the positions $k\geq i$ below or on the diagonal of the matrix.

\begin{example} \label{Ex1}
{\rm Let us look at an example of the kind that occurs in the proof of Proposition \ref{maxd=3mindd=1}.
The entries $a_j$ and $a_j'$ of 
\begin{align*}
a&=(0,3,5,6,8,10,12,15,16,20,23),\;\; 
a'=(0,3,5,7,9,11,13,15,16,20,23) \in A(24,11),
\end{align*}
pairwise coincide for $1\leq j\leq 3=:u$ and for $v:=8 \leq j \leq 11=:r$.
The delta tableaux $\mathfrak{A}$ and $\mathfrak{A}'$ of $a$ and $a'$, respectively, are
\begin{center}
\begin{tabular}{cccccccccc}
\cellcolor[gray]{0.7}3&\cellcolor[gray]{0.7} 5&\cellcolor[gray]{0.9}6&\cellcolor[gray]{0.9}8&\cellcolor[gray]{0.9}10&\cellcolor[gray]{0.9}12&\cellcolor[gray]{0.7}15&\cellcolor[gray]{0.7}16&\cellcolor[gray]{0.7}20&\cellcolor[gray]{0.7}23 \\
&\cellcolor[gray]{0.7} 2&\cellcolor[gray]{0.9}3&\cellcolor[gray]{0.9}5&\cellcolor[gray]{0.9}7&\cellcolor[gray]{0.9}9&\cellcolor[gray]{0.7}12&\cellcolor[gray]{0.7}13&\cellcolor[gray]{0.7}17&\cellcolor[gray]{0.7}20 \\
&&\cellcolor[gray]{0.9}1&\cellcolor[gray]{0.9}3&\cellcolor[gray]{0.9}5&\cellcolor[gray]{0.9}7&\cellcolor[gray]{0.7}10&\cellcolor[gray]{0.7}11&\cellcolor[gray]{0.7}15&\cellcolor[gray]{0.7}18\\ 
&\multicolumn{1}{c}{}&&\cellcolor[gray]{0.7}2&\cellcolor[gray]{0.7}4&\cellcolor[gray]{0.7}6&\cellcolor[gray]{0.9}9&\cellcolor[gray]{0.9}10&\cellcolor[gray]{0.9}14&\cellcolor[gray]{0.9}17\\
&\multicolumn{1}{c}{}&&&\cellcolor[gray]{0.7}2&\cellcolor[gray]{0.7}4&\cellcolor[gray]{0.9}7&\cellcolor[gray]{0.9}8&\cellcolor[gray]{0.9}12&\cellcolor[gray]{0.9}15\\
&\multicolumn{1}{c}{}&&&&\cellcolor[gray]{0.7}2&\cellcolor[gray]{0.9}5&\cellcolor[gray]{0.9}6&\cellcolor[gray]{0.9}10&\cellcolor[gray]{0.9}13\\
&\multicolumn{1}{c}{}&&&&&\cellcolor[gray]{0.9}3&\cellcolor[gray]{0.9}4&\cellcolor[gray]{0.9}8&\cellcolor[gray]{0.9}11\\ 
&\multicolumn{1}{c}{}&&&&\multicolumn{1}{c}{}&&\cellcolor[gray]{0.7}1&\cellcolor[gray]{0.7}5&\cellcolor[gray]{0.7}8\\
&\multicolumn{1}{c}{}&&&&\multicolumn{1}{c}{}&&&\cellcolor[gray]{0.7}4&\cellcolor[gray]{0.7}7\\
&\multicolumn{1}{c}{}&&&&\multicolumn{1}{c}{}&&&&\cellcolor[gray]{0.7}3\\
\end{tabular}
\quad\quad
\begin{tabular}{cccccccccc}
\cellcolor[gray]{0.7}3&\cellcolor[gray]{0.7}5&\cellcolor[gray]{0.9}7&\cellcolor[gray]{0.9}9&\cellcolor[gray]{0.9}11&\cellcolor[gray]{0.9}13&\cellcolor[gray]{0.7}15&\cellcolor[gray]{0.7}16&\cellcolor[gray]{0.7}20&\cellcolor[gray]{0.7}23 \\
&\cellcolor[gray]{0.7}2&\cellcolor[gray]{0.9}4&\cellcolor[gray]{0.9}6&\cellcolor[gray]{0.9}8&\cellcolor[gray]{0.9}10&\cellcolor[gray]{0.7}12&\cellcolor[gray]{0.7}13&\cellcolor[gray]{0.7}17&\cellcolor[gray]{0.7}20 \\
&&\cellcolor[gray]{0.9}2&\cellcolor[gray]{0.9}4&\cellcolor[gray]{0.9}6&\cellcolor[gray]{0.9}8&\cellcolor[gray]{0.7}10&\cellcolor[gray]{0.7}11&\cellcolor[gray]{0.7}15&\cellcolor[gray]{0.7}18\\ 
&\multicolumn{1}{c}{}&&\cellcolor[gray]{0.7}2&\cellcolor[gray]{0.7}4&\cellcolor[gray]{0.7}6&\cellcolor[gray]{0.9}8&\cellcolor[gray]{0.9}9&\cellcolor[gray]{0.9}13&\cellcolor[gray]{0.9}16\\
&\multicolumn{1}{c}{}&&&\cellcolor[gray]{0.7}2&\cellcolor[gray]{0.7}4&\cellcolor[gray]{0.9}6&\cellcolor[gray]{0.9}7&\cellcolor[gray]{0.9}11&\cellcolor[gray]{0.9}14\\
&\multicolumn{1}{c}{}&&&&\cellcolor[gray]{0.7}2&\cellcolor[gray]{0.9}4&\cellcolor[gray]{0.9}5&\cellcolor[gray]{0.9}9&\cellcolor[gray]{0.9}12\\
&\multicolumn{1}{c}{}&&&&&\cellcolor[gray]{0.9}2&\cellcolor[gray]{0.9}3&\cellcolor[gray]{0.9}7&\cellcolor[gray]{0.9}10\\ 
&\multicolumn{1}{c}{}&&&&\multicolumn{1}{c}{}&&\cellcolor[gray]{0.7}1&\cellcolor[gray]{0.7}5&\cellcolor[gray]{0.7}8\\
&\multicolumn{1}{c}{}&&&&\multicolumn{1}{c}{}&&&\cellcolor[gray]{0.7}4&\cellcolor[gray]{0.7}7\\
&\multicolumn{1}{c}{}&&&&\multicolumn{1}{c}{}&&&&\cellcolor[gray]{0.7}3\\
\end{tabular}
\end{center}
If we wish to evaluate $h_{p,11}(a)-h_{p,11}(a')$, we have to subtract summands of type $p^{-a_{k,i}'}$
from summands of type $p^{-a_{k,i}}$. 
We observe that the entries inside the dark-coloured areas (three triangles and the upper right rectangle) of $\mathfrak{A}$ one by one coincide with those of $\mathfrak{A}'$, which means that the corresponding summands annihilate each other. Consequently, we only have to consider the entries in each of the two remaining light-coloured rectangles of $\mathfrak{A}$ and $\mathfrak{A}'$, where differences $a_{k,i}-a'_{k,i}$ are $-1$ between the upper left rectangles and $+1$ between the lower right rectangles. Note that the three vertical blocks in the tableaux have column indices $i$ running in the intervals $2\leq i \leq u$ and $u+1\leq i \leq v-1$ and $v\leq i \leq r$, while the row indices of the three horizontal blocks have ranges $1\leq k \leq u$ and $u+1\leq k \leq v-1$ and $v \leq k \leq r-1$.
}
\end{example}

The following result is a corollary to Proposition 3.1 in \cite{SA3}. Yet, in view of the new perspective of this exposition and for the convenience of the reader, we provide a short proof of it.
\begin{proposition} \label{maxd=3mindd=1}
Let $p$ be a prime. 
For some integers $2\leq r \leq s$ and $1\leq u < v \leq r-1$, let $d=(d_1,\ldots,d_{r-1})\in D(s,r)$ have the following properties:
\begin{itemize}
\item[(i)] $(d_u,d_v) \in \{(1,3),(3,1)\}$,
\item[(ii)]  $d_j=2$ for $u<j<v$.
\end{itemize}
Then $d'=(d_1',\ldots,d_{r-1}')$, defined by
\[
d_j' := \left\{ \begin{array}{ll}
                   2 &  \mbox{ for $j=u$ and $j=v$,} \\
                   d_j &  \mbox{ otherwise,} 
                   \end{array}  
                   \right.    
                   \]
satisfies $d'\in D(s,r)$ and $\Ene{{\cal D}(\delta_{r}^{-1}(d'))}{p^s} > \Ene{{\cal D}(\delta_r^{-1}(d))}{p^s}$.
\end{proposition}

{\sc Proof. }
By condition (i) and the symmetry property
(\ref{symmhp}) of $h_{p,r}$ we may assume w.l.o.g. that $d_u=1$ and $d_v=3$, i.e.
\begin{align*}
d &= (d_1,\ldots,d_{u-1},1,2,\ldots,2,3,d_{v+1},\ldots,d_{r-1}) \in D(s,r),\\[-.2in]
\intertext{and}\\[-.45in]  
d' &= (d_1,\ldots,d_{u-1},2,2,\ldots,2,2,d_{v+1},\ldots,d_{r-1}).
\end{align*}
Clearly, $\sum d_j' = \sum d_j =s-1$, hence $d'\in D(s,r)$.
For $a=(a_1,\ldots,a_r) := \delta_r^{-1}(d) \in A(s,r)$, this means that
$a' = (a_1',\ldots, a_r') := \delta_{r}^{-1}(d') \in A(s,r)$ satisfies
\begin{equation} \label{aprime}
a_j' = \left\{ \begin{array}{ll}
                   a_j &  \mbox{ for $1\leq j \leq u$,} \\
                   a_j+1 &  \mbox{ for $u+1 \leq j \leq v$,} \\
                   a_j &  \mbox{ for $v+1 \leq j \leq r$.} 
                   \end{array}  
                   \right.                  
\end{equation}
By definition of the function $h_{p,r}$ in (\ref{basishp}), we have
\[
h_{p,r}(a) - h_{p,r}(a') = 
 \sum_{k=1}^{r-1}\sum_{i=k+1}^{r} \,\left(\frac{1}{p^{a_i-a_k}}  -   \frac{1}{p^{a'_i-a'_k}}\right). 
\]
By use of (\ref{aprime}), comparison of the delta tableaux $\mathfrak{A}=(a_{k,i})_{1\leq k<i\leq r}$ and $\mathfrak{A}'=(a'_{k,i})_{1\leq k<i\leq r}$ (cf. Example \ref{Ex1}), where $a_{k,i}:= a_i-a_k$ and $a'_{k,i}:= a'_i-a'_k$ for $1\leq k<i\leq r$, reveals that 
\begin{align*} 
\begin{split}
h_{p,r}(a) - h_{p,r}(a') 
&= \sum_{k=1}^u \sum_{i=u+1}^v \,\left(\frac{1}{p^{a_{k,i}}}  -   \frac{1}{p^{a'_{k,i}}}\right)
 +  
\sum_{k=u+1}^v \sum_{i=v+1}^r \,\left(\frac{1}{p^{a_{k,i}}}  -   \frac{1}{p^{a'_{k,i}}}\right) \\
&= \sum_{k=1}^u \sum_{i=u+1}^v \,\left(\frac{1}{p^{a_{k,i}}}  -   \frac{1}{p^{a_{k,i}+1}}\right)
 +  
\sum_{k=u+1}^v \sum_{i=v+1}^r \,\left(\frac{1}{p^{a_{k,i}}}  -   \frac{1}{p^{a_{k,i}-1}}\right) \\
&= \left(1-\frac{1}{p} \right) \sum_{k=1}^u \sum_{i=u+1}^v \,\frac{1}{p^{a_{k,i}}}
 +  
(1-p) \sum_{k=u+1}^v \sum_{i=v+1}^r \,\frac{1}{p^{a_{k,i}}}.  
\end{split}
\end{align*}
Since assumption (ii) implies that
\begin{equation}  \label{diffaj}
a_j = a_{u+1}+2(j-(u+1)) \quad\quad\quad (u+1 \leq j \leq v),
\end{equation}
it follows that 
\begin{align} \label{diffhp2}
\begin{split}
h_{p,r}(a) - h_{p,r}(a') 
&= \left(1-\frac{1}{p} \right) \sum_{k=1}^u p^{a_k} \sum_{i=u+1}^v \,\frac{1}{p^{a_{u+1}+2(i-(u+1))}}  \\
 & \quad\quad\quad\quad +  
(1-p) \sum_{k=u+1}^v p^{a_{u+1}+2(k-(u+1))}\sum_{i=v+1}^r \,\frac{1}{p^{a_i}}  \\
&= \left( 1- \frac{1}{p}\right)  \frac{1}{p^{a_{u+1}}} \sum_{k=1}^u p^{a_k} \sum_{i=0}^{v-u-1} \,\frac{1}{p^{2i}}  \\
 & \quad\quad\quad\quad -  
(p-1) p^{a_{u+1}} \sum_{k=0}^{v-u-1} p^{2k} \sum_{i=v+1}^r \,\frac{1}{p^{a_i}}.
\end{split}
\end{align}
The identity $a_v = a_{u+1} +2(v-(u+1))$, being a special case of (\ref{diffaj}), yields
\[  \frac{1}{p^{a_{u+1}}} \sum_{i=0}^{v-u-1} \,\frac{1}{p^{2i}} = 
  \frac{1}{p^{a_{u+1}}}\cdot \frac{1}{p^{2(v-u-1)}} \sum_{i=0}^{v-u-1} \,p^{2i} =
  \frac{1}{p^{a_v}} \sum_{i=0}^{v-u-1} \,p^{2i}\, .
\]
Inserting this into (\ref{diffhp2}), we obtain
\begin{align} \label{diffhp3}
\begin{split}
h_{p,r}(a) - h_{p,r}(a') 
&= \left(\frac{1}{p^{a_v+1}} \sum_{k=1}^u p^{a_k}  -  
 p^{a_{u+1}} \sum_{i=v+1}^r \,\frac{1}{p^{a_i}}
 \right) (p-1)\sum_{i=0}^{v-u-1} \,p^{2i} \\
&=\left(\frac{1}{p^{a_v+1}} \sum_{k=1}^u p^{a_k}  -  
 p^{a_{u+1}} \sum_{i=v+1}^r \,\frac{1}{p^{a_i}}
 \right) \frac{p^{2(v-u)}-1}{p+1} \, .
\end{split}
\end{align}
Since $d_u=1$ and $d_v=3$ by assumption, we have $a_{u+1}=a_u+1$ and 
$a_{v+1}=a_v+3$. 
Therefore, $\sum_{k=1}^u p^{a_k} \geq p^{a_u} = p^{a_{u+1}-1}$ and 
\[
\sum_{i=v+1}^r \frac{1}{p^{a_i}}
< \frac{1}{p^{a_{v+1}}} \sum_{i=0}^{\infty} \frac{1}{p^i} = \frac{p}{(p-1)p^{a_{v+1}}} =  \frac{1}{(p-1)p^{a_v +2}}\, .
\]
With this, (\ref{diffhp3}) yields
\[
h_{p,r}(a) - h_{p,r}(a') 
> \left(p^{a_{u+1}-a_v-2}  -  
 \frac{p^{a_{u+1}-a_v-2}}{p-1} \right) \frac{p^{2(v-u)}-1}{p+1}
= \frac{p-2}{p^2-1}\cdot p^{a_{u+1}-a_v-2} (p^{2(v-u)}-1)
 \, , 
\]
hence $h_{p,r}(a) - h_{p,r}(a') > 0$ for each prime $p$. By Lemma \ref{hilfssatz}, we finally obtain the desired conclusion.

\qed

In Proposition \ref{maxd=3mindd=1}, tableaux of the same size were ``subtracted'' from each other. Now we study the case where a given tableau has to be compared with a smaller one.
For that purpose, we introduce the $(u,v)$-\textit{derivative} of an admissible tuple.
\begin{definition}
Let $a\in A(s,r)$ for some $3\leq r \leq s$, and let $1\leq u < v \leq r-1$ be arbitrary integers. Then $\partial_{u,v}(a)=(a'_1,\ldots,a'_{r-1}) \in A(s,r-1)$, defined by setting
\begin{equation*}   \label{aprimeTRL}
a_j' := \left\{ \begin{array}{ll}
                   a_j &  \mbox{ for $1\leq j \leq u$,} \\
                   a_j+1 &  \mbox{ for $u+1 \leq j \leq v-1$,} \\
                   a_{j+1} &  \mbox{ for $v \leq j \leq r-1$,} 
                   \end{array}  
                   \right.                  
\end{equation*}
is called the $(u,v)$-\textit{derivative} of $a$.
\end{definition}

\begin{lemma}[Tableau Reduction Lemma]{\label{tablred}}
Let $p$ be a prime, and let $a=(a_1,\ldots,a_r) \in A(s,r)$ be admissible for some $3\leq r \leq s$. 
Moreover, let $1\leq u < v \leq r-1$ be arbitrary integers.
If $a_{j+1}-a_j=2$ for $u+1\leq j \leq v-1$, then 
\begin{align}  \label{deltahfinal}
\begin{split}
h_{p,r}(a) - h_{p,r-1}(\partial_{u,v}(a)) &= \frac{1}{p+1}\left(p+\frac{1}{p^{2(v-u-1)}}\right) \left(\frac{U}{p^{a_{u+1}}}+p^{a_v} V\right) \\
&\quad + \frac{1}{p^2-1}\left(1- \frac{1}{p^{2(v-u-1)}}\right),
\end{split}
\end{align}
where
\[
U := \sum_{k=1}^u p^{a_k}  \quad \mbox{ and } \quad
V:= \sum_{i=v+1}^r \frac{1}{p^{a_i}}.
\]
\end{lemma}

{\sc Proof. }
%
By definition of the functions $h_{p,r}$ and $h_{p,r-1}$ in 
(\ref{basishp}), 
we have for $a' = (a_1',\ldots, a'_{r-1}):=\partial_{u,v}(a)$
\[
\Delta h:= h_{p,r}(a) - h_{p,r-1}(\partial_{u,v}(a)) = 
\sum_{k=1}^{r-1}\sum_{i=k+1}^{r} \,\frac{1}{p^{a_i-a_k}} -   \sum_{k=1}^{r-2}\sum_{i=k+1}^{r-1} \, \frac{1}{p^{a'_i-a'_k}}.
\]
Using (\ref{aprimeTRL}), 
we obtain the following two delta tableaux $\mathfrak{A}=(a_{k,i})_{1\leq k<i\leq r}$ and $\mathfrak{A}'=(a'_{k,i})_{1\leq k<i\leq r-1}$, where $a_{k,i}:= a_i-a_k$  ($1\leq k<i\leq r$) and $a'_{k,i}:= a'_i-a'_k$ for $1\leq k<i\leq r-1$:
\footnotesize
\bigskip
\medskip
\newline
\begin{tabular}{cccccccccccc}
\cellcolor[gray]{0.7}$a_{1,2}$&\cellcolor[gray]{0.7}$\ldots\ldots$&\cellcolor[gray]{0.7}$a_{1,u}$&\cellcolor[gray]{0.9}$a_{1,u+1}$ &\multicolumn{3}{c}{\cellcolor[gray]{0.9}$\ldots\ldots\ldots\ldots\ldots\ldots\ldots\ldots$}
&\cellcolor[gray]{0.9}$a_{1,v}$&\cellcolor[gray]{0.7}$a_{1,v+1}$&\multicolumn{2}{c}{\cellcolor[gray]{0.7}$\ldots\ldots\ldots\ldots\ldots$}&\cellcolor[gray]{0.7}$a_{1,r}$ \\
&\cellcolor[gray]{0.7}$\ddots$&\cellcolor[gray]{0.7}$\vdots$&\cellcolor[gray]{0.9}\vdots&\cellcolor[gray]{0.9}&\cellcolor[gray]{0.9}&\cellcolor[gray]{0.9}&\cellcolor[gray]{0.9}$\vdots$&\cellcolor[gray]{0.7}$\vdots$&\cellcolor[gray]{0.7}&\cellcolor[gray]{0.7}&\cellcolor[gray]{0.7}$\vdots$ \\
&&\cellcolor[gray]{0.7}$a_{u-1,u}$&\cellcolor[gray]{0.9}&\cellcolor[gray]{0.9}&\cellcolor[gray]{0.9}&\cellcolor[gray]{0.9}&\cellcolor[gray]{0.9}&\cellcolor[gray]{0.7}&\cellcolor[gray]{0.7}& \cellcolor[gray]{0.7}&\cellcolor[gray]{0.7}\\ 
&&&\cellcolor[gray]{0.9}$a_{u,u+1}$&\multicolumn{3}{c}{\cellcolor[gray]{0.9}$\ldots\ldots\ldots\ldots\ldots\ldots\ldots\ldots$}&\cellcolor[gray]{0.9}$a_{u,v}$ &\cellcolor[gray]{0.7}$a_{u,v+1}$&\multicolumn{2}{c}{\cellcolor[gray]{0.7}$\ldots\ldots\ldots\ldots\ldots$}&\cellcolor[gray]{0.7}$a_{u,r}$\\  
\multicolumn{4}{c}{}&\cellcolor[gray]{0.7}$a_{u+1,u+2}$&{\cellcolor[gray]{0.7}$\ldots$}&{\cellcolor[gray]{0.7}$a_{u+1,v-1}$}&$a_{u+1,v}$&\cellcolor[gray]{0.9}$a_{u+1,v+1}$&\multicolumn{2}{c}{\cellcolor[gray]{0.9}$\ldots\ldots\ldots\ldots\ldots$}&\cellcolor[gray]{0.9}$a_{u+1,r}$ \\ 
\multicolumn{4}{c}{}&&\cellcolor[gray]{0.7}$\ddots$&\cellcolor[gray]{0.7}$\vdots$&$\vdots$&\cellcolor[gray]{0.9}&\cellcolor[gray]{0.9}&\cellcolor[gray]{0.9}&\cellcolor[gray]{0.9}\\[-.075in]
\multicolumn{4}{c}{}&&&\cellcolor[gray]{0.7}$a_{v-2,v-1}$&$a_{v-2,v}$&\cellcolor[gray]{0.9}$\vdots$&\cellcolor[gray]{0.9}&\cellcolor[gray]{0.9}&\cellcolor[gray]{0.9}$\vdots$\\
\multicolumn{7}{c|}{}&$a_{v-1,v}$&\cellcolor[gray]{0.9}&\cellcolor[gray]{0.9} &\cellcolor[gray]{0.9}&\cellcolor[gray]{0.9}\\ \cline{8-8}
\multicolumn{4}{c}{}&&&&&\cellcolor[gray]{0.9}$a_{v,v+1}$&\multicolumn{2}{c}{\cellcolor[gray]{0.9}$\ldots\ldots\ldots\ldots\ldots$}&\cellcolor[gray]{0.9}$a_{v,r}$ \\ 
\multicolumn{9}{c}{}&\cellcolor[gray]{0.7}$a_{v+1,v+2}$&\cellcolor[gray]{0.7}$\ldots$&\cellcolor[gray]{0.7}$a_{v+1,r}$\\
\multicolumn{9}{c}{}&&\cellcolor[gray]{0.7}$\ddots$&\cellcolor[gray]{0.7}$\vdots$\\
\multicolumn{9}{c}{}&&&\cellcolor[gray]{0.7}$a_{r-1,r}$
\end{tabular}
\bigskip
\bigskip
\medskip
\newline
\begin{tabular}{ccccccccccc}
\cellcolor[gray]{0.7}$a_{1,2}$&\cellcolor[gray]{0.7}$\ldots\ldots$&\cellcolor[gray]{0.7}
$a_{1,u}$&\cellcolor[gray]{0.9}$a'_{1,u+1}$&\multicolumn{2}{c}{\cellcolor[gray]{0.9}$\ldots\ldots\ldots\ldots$}&\cellcolor[gray]{0.9}$a'_{1,v-1}$&\cellcolor[gray]{0.7}
$a_{1,v+1}$&\multicolumn{2}{c}{\cellcolor[gray]{0.7}$\ldots\ldots\ldots\ldots\ldots$}&\cellcolor[gray]{0.7}$a_{1,r}$ \\
&\cellcolor[gray]{0.7}$\ddots$&\cellcolor[gray]{0.7}$\vdots$&\cellcolor[gray]{0.9}$\vdots$&\cellcolor[gray]{0.9}&\cellcolor[gray]{0.9}&\cellcolor[gray]{0.9}$\vdots$&\cellcolor[gray]{0.7}$\vdots$&\cellcolor[gray]{0.7}&\cellcolor[gray]{0.7}&\cellcolor[gray]{0.7}$\vdots$ \\
&&\cellcolor[gray]{0.7}$a_{u-1,u}$&\cellcolor[gray]{0.9}&\cellcolor[gray]{0.9}&\cellcolor[gray]{0.9}&\cellcolor[gray]{0.9}&\cellcolor[gray]{0.7}&\cellcolor[gray]{0.7}&\cellcolor[gray]{0.7} &\cellcolor[gray]{0.7}\\ 
&&&\cellcolor[gray]{0.9}$a'_{u,u+1}$&\multicolumn{2}{c}{\cellcolor[gray]{0.9}$\ldots\ldots\ldots\ldots$}&\cellcolor[gray]{0.9}$a'_{u,v-1}$&\cellcolor[gray]{0.7}$a_{u,v+1}$&\multicolumn{2}{c}{\cellcolor[gray]{0.7}$\ldots\ldots\ldots\ldots\ldots$}&\cellcolor[gray]{0.7}$a_{u,r}$\\ 
\multicolumn{4}{c}{}&\cellcolor[gray]{0.7}$a_{u+1,u+2}$&\cellcolor[gray]{0.7}$\ldots$\cellcolor[gray]{0.7}&\cellcolor[gray]{0.7}$a_{u+1,v-1}$&\cellcolor[gray]{0.9}$a'_{u+1,v}$&\multicolumn{2}{c}{\cellcolor[gray]{0.9}\ldots\ldots\ldots\ldots\ldots}&\cellcolor[gray]{0.9}$a'_{u+1,r-1}$ \\
\multicolumn{4}{c}{}&&\cellcolor[gray]{0.7}$\ddots$&\cellcolor[gray]{0.7}$\vdots$&\cellcolor[gray]{0.9}$\vdots$&\cellcolor[gray]{0.9}&\cellcolor[gray]{0.9}&\cellcolor[gray]{0.9}$\vdots$\\
\multicolumn{4}{c}{}&&&\cellcolor[gray]{0.7}$a_{v-2,v-1}$&\cellcolor[gray]{0.9}$a'_{v-2,v}$&\multicolumn{2}{c}{\cellcolor[gray]{0.9}$\ldots\ldots\ldots\ldots\ldots$}&\cellcolor[gray]{0.9}$a'_{v-2,r-1}$\\
\multicolumn{3}{c}{}&&&&&\cellcolor[gray]{0.9}$a'_{v-1,v}$&\multicolumn{2}{c}{\cellcolor[gray]{0.9}$\ldots\ldots\ldots\ldots\ldots$}&\cellcolor[gray]{0.9}$a'_{v-1,r}$\\ 
\multicolumn{8}{c}{}&\cellcolor[gray]{0.7}$a_{v+1,v+2}$&\cellcolor[gray]{0.7}$\ldots$&\cellcolor[gray]{0.7}$a_{v+1,r}$\\
\multicolumn{8}{c}{}&&\cellcolor[gray]{0.7}$\ddots$&\cellcolor[gray]{0.7}$\vdots$\\
\multicolumn{8}{c}{}&&&\cellcolor[gray]{0.7}$a_{r-1,r}$
\end{tabular}
\normalsize
\bigskip
\bigskip\newline
Observe that the entries inside the dark-coloured areas (three triangles and the upper right rectangle) of $\mathfrak{A}$ one by one equal the corresponding entries of $\mathfrak{A}'$.
For that reason, we only have to consider the entries in each of the two remaining light-coloured rectangles of $\mathfrak{A}$ and $\mathfrak{A}'$ and finally the entries in the white-coloured column inside $\mathfrak{A}$.
Hence
\begin{align*}
\Delta h&= 
\left( \sum_{k=1}^{u}\sum_{i=u+1}^{v} 
+ \sum_{k=u+1}^{v}\sum_{i=v+1}^{r} \right) \frac{1}{p^{a_{k,i}}} + \sum_{k=u+1}^{v-1}\, \frac{1}{p^{a_{k,v}}} \; - 
\left( \sum_{k=1}^{u}\sum_{i=u+1}^{v-1} + \sum_{k=u+1}^{v-1}\sum_{i=v}^{r-1}\right) \frac{1}{p^{a'_{k,i}}}.
 \end{align*}
Then (\ref{aprimeTRL}) 
yields
\begin{align} \label{diffhp7}
\begin{split}
\Delta h
&=  \sum_{k=1}^{u}\left(\sum_{i=u+1}^{v-1} \, \left(\frac{1}{p^{a_{k,i}}}  -\frac{1}{p^{a'_{k,i}}}\right) + \frac{1}{p^{a_{k,v}}}\right)\\ 
&\quad\quad 
+ \sum_{k=u+1}^{v-1}\left(\sum_{i=v+1}^{r} \frac{1}{p^{a_{k,i}}}  -\sum_{i=v}^{r-1} \frac{1}{p^{a'_{k,i}}}\right) + \sum_{i=v+1}^{r}\,\frac{1}{p^{a_{v,i}}}  +
\sum_{k=u+1}^{v-1}\, \frac{1}{p^{a_{k,v}}}       \\
&=  \sum_{k=1}^{u}\left(\sum_{i=u+1}^{v-1} \, \left(\frac{1}{p^{a_{k,i}}}  -\frac{1}{p^{a_{k,i}+1}}\right)\right) + \sum_{k=1}^u\frac{1}{p^{a_{k,v}}}\\ 
&\quad\quad 
+ \sum_{k=u+1}^{v-1}\left(\sum_{i=v+1}^{r} \frac{1}{p^{a_{k,i}}}  -\sum_{i=v+1}^{r} \frac{1}{p^{a_{k,i}-1}}\right) + \sum_{i=v+1}^{r}\,\frac{1}{p^{a_{v,i}}} +
\sum_{k=u+1}^{v-1}\, \frac{1}{p^{a_{k,v}}}      \\
&= \frac{p-1}{p} \sum_{k=1}^{u}p^{a_k} \sum_{i=u+1}^{v-1} \, \frac{1}{p^{a_i}}  + \frac{1}{p^{a_v}}\sum_{k=1}^u p^{a_k}\\ 
&\quad\quad 
+ (1-p) \sum_{k=u+1}^{v-1} p^{a_k} \sum_{i=v+1}^{r} \frac{1}{p^{a_{i}}}  
+ p^{a_v} \sum_{i=v+1}^{r}\,\frac{1}{p^{a_i}}+
\sum_{k=u+1}^{v-1}\, \frac{1}{p^{a_{k,v}}}     \\
&= 
\frac{p-1}{p}\, U \sum_{i=u+1}^{v-1} \, \frac{1}{p^{a_i}}  + \frac{1}{p^{a_v}}\, U
- (p-1) V \sum_{k=u+1}^{v-1} p^{a_k} + p^{a_v} V  +\sum_{k=u+1}^{v-1}\, \frac{1}{p^{a_{k,v}}}.
\end{split}
\end{align}

Our assumption 
$a_{j+1}-a_j = 2$ for $u+1\leq j \leq v-1$ implies that 
\begin{equation}  \label{diffaj7}
a_j = a_{u+1}+2(j-(u+1)) \quad\quad\quad (u+1 \leq j \leq v).
\end{equation}
It follows that
\begin{align*}
\sum_{k=u+1}^{v-1} \,p^{a_k} &= p^{a_{v-1}} \sum_{k=u+1}^{v-1} \,\frac{1}{p^{a_{v-1}-a_k}} = p^{a_{v-1}} \sum_{k=0}^{v-u-2} \frac{1}{p^{2k}} = p^{a_{v-1}} \frac{p^2}{p^2-1}\left(1- \frac{1}{p^{2(v-u-1)}}\right),  \\
\sum_{i=u+1}^{v-1} \,\frac{1}{p^{a_i}} &= \frac{1}{p^{a_{u+1}}} \sum_{k=u+1}^{v-1} \,p^{a_i-a_{u+1}} = \frac{1}{p^{a_{u+1}}} \sum_{i=0}^{v-u-2} \frac{1}{p^{2i}} = \frac{1}{p^{a_{u+1}}}\cdot \frac{p^2}{p^2-1}\left(1- \frac{1}{p^{2(v-u-1)}}\right),
\\[-.15in]
\intertext{and}\\[-.45in] 
\sum_{i=u+1}^{v-1}\, \frac{1}{p^{a_{k,v}}} &= \sum_{i=u+1}^{v-1}\, \frac{1}{p^{2(v-k)}} = \sum_{k=1}^{v-u-1} \frac{1}{p^{2k}} = \frac{1}{p^2-1}\left(1- \frac{1}{p^{2(v-u-1)}}\right).
\end{align*}
If we insert these identities into (\ref{diffhp7}) and use 
$a_v=a_{u+1}+2(v-u-1)$, being a special case of 
(\ref{diffaj7}), we obtain
\begin{align} \label{lastident}
\begin{split}
\Delta h
&=
\left( \frac{1}{p^{a_v}} + \frac{1}{p^{a_{u+1}}}\cdot \frac{p}{p+1}\left(1- \frac{1}{p^{2(v-u-1)}}\right)\right)U \\
&\quad\quad\quad+ \left( p^{a_v} - p^{a_{v-1}}\,\frac{p^2}{p+1}\left(1- \frac{1}{p^{2(v-u-1)}}\right)\right)V + \frac{1}{p^2-1} \left(1- \frac{1}{p^{2(v-u-1)}}\right) \\
&=
\left( \frac{1}{p^{2(v-u-1)}} + \frac{p}{p+1}\left(1- \frac{1}{p^{2(v-u-1)}}\right)\right)\frac{U}{p^{a_{u+1}}} \\
&\quad\quad\quad+ \left(1 - \frac{p^{2-(a_v-a_{v-1})}}{p+1}\left(1- \frac{1}{p^{2(v-u-1)}}\right)\right)p^{a_v} V
+\frac{1}{p^2-1} \left(1- \frac{1}{p^{2(v-u-1)}}\right). 
\end{split}
\end{align}
In case $v\geq u+2$ we know by (\ref{diffaj7}) that $a_v-a_{v-1}=2$, and (\ref{deltahfinal}) follows.
To complete the proof we observe that identity (\ref{deltahfinal}) is a trivial consequence of (\ref{lastident}) in case $v=u+1$. 

\qed

\begin{proposition}{\label{twice1}}
Let $p$ be a prime. 
For some integers $3\leq r \leq s$ and $1\leq u < v \leq r-1$, let $d=(d_1,\ldots,d_{r-1})\in D(s,r)$ have the following properties:
\begin{itemize}
\item[(i)] $d_u=d_v=1$, 
\item[(ii)]  $d_j=2$ for $u<j<v$.
\end{itemize}
Then $d'=(d_1',\ldots,d_{r-2}')$, defined by
\[
d_j' := \left\{ \begin{array}{ll}
                   d_j &  \mbox{ for $1\leq j \leq u-1$,} \\
                   2 &  \mbox{ for $u\leq j \leq v-1$,} \\
                   d_{j+1} &  \mbox{ for $v\leq j \leq r-2$,} 
                   \end{array}  
                   \right.    
\]
satisfies $d'\in D(s,r-1)$ and $\Ene{{\cal D}(\delta_{r-1}^{-1}(d'))}{p^s} > \Ene{{\cal D}(\delta_r^{-1}(d))}{p^s}$
with precisely one exceptional case, namely
$p=2$, $u=1$ and $v=r-1$, i.e. $p=2$ and $d=(1,2,2,\ldots,2,1)$, where $d'=(2,2,\ldots,2,2)\in D(s,r-1)$ and
$\Ene{{\cal D}(\delta_{r-1}^{-1}(d'))}{2^s} = \Ene{{\cal D}(\delta_r^{-1}(d))}{2^s}$.
\end{proposition}

{\sc Proof. } Since $\sum d'_j = \sum d_j =s-1$, we have $d'\in D(s,r-1)$.
We set $a=(a_1,\ldots,a_r) := \delta_r^{-1}(d) \in A(s,r)$. Then it is easy to check that  
$a':=\delta_{r-1}^{-1}(d') \in A(s,r-1)$  is the $(u,v)$-derivative of $a$. 
By virtue of condition (ii), we can apply the Tableau Reduction Lemma \ref{tablred}. Before evaluating (\ref{deltahfinal}), we observe that 
\[
U= \sum_{k=1}^u p^{a_k} \geq p^{a_u} \quad \mbox{ and } \quad
V= \sum_{i=v+1}^r \frac{1}{p^{a_i}} \geq \frac{1}{p^{a_{v+1}}},
\]
where equality in both cases simultaneously holds if and only if $u=1$ and $v=r-1$. 
Inserting this into (\ref{deltahfinal}) and making use of (i), that is $a_{u+1}-a_u=d_u=1$ and $a_{v+1}-a_v=d_v=1$  by definition of $\delta_r$, we obtain
\begin{align*}  
\begin{split}
h_{p,r}(a) - h_{p,r-1}(a') &= h_{p,r}(a) - h_{p,r-1}(a')  \\
&\geq \frac{1}{p+1}\left(p+\frac{1}{p^{2(v-u-1)}}\right) \left(\frac{p^{a_u}}{p^{a_{u+1}}}+\frac{p^{a_v}}{p^{a_{v+1}}}\right) + \frac{1}{p^2-1}\left(1- \frac{1}{p^{2(v-u-1)}}\right)\\
&= \frac{2}{p(p+1)}\left(p+\frac{1}{p^{2(v-u-1)}}\right) + \frac{1}{p^2-1}\left(1- \frac{1}{p^{2(v-u-1)}}\right) \\
&= \frac{2p-1}{p^2-1} + \frac{p-2}{p^{2v-2u-1}(p^2-1)}\\
&\geq \frac{1}{p-1}
\end{split}
\end{align*}
for all primes $p$, but equality in the final step only for $p=2$. Altogether, we
have shown that $h_{p,r}(a) - h_{p,r-1}(a') > \frac{1}{p-1}$ in all cases with the unique exception excluded in the Proposition, where equality holds.
The proof is completed by  Lemma \ref{hilfssatz}. 
%

\qed

\begin{proposition}{\label{twice3}}
Let $p$ be a prime. 
For some integers $4\leq r \leq s$ and $1\leq u < v \leq r-2$, let $d'=(d'_1,\ldots,d'_{r-2})\in D(s,r-1)$ have the following properties:
\begin{itemize}
\item[(i)] $d'_u=d'_v=3$, 
\item[(ii)]  $d'_j=2$ for $u<j<v$.
\end{itemize}
Then $d=(d_1,\ldots,d_{r-1})$, defined by
\[
d_j := \left\{ \begin{array}{ll}
                   d'_j &  \mbox{ for $1\leq j \leq u-1$ or $u+1\leq j \leq v-1$,} \\
                   2 &  \mbox{ for $u \leq j \leq v+1$,} \\
                   d'_{j-1} &  \mbox{ for $v+2 \leq j \leq r-1$,} 
                   \end{array}  
                   \right.    
\]
satisfies $d\in D(s,r)$ and $\Ene{{\cal D}(\delta_{r}^{-1}(d))}{p^s} > \Ene{{\cal D}(\delta_{r-1}^{-1}(d'))}{p^s}$.
\end{proposition}

{\sc Proof. }
Since $\sum d_j = \sum d'_j =s-1$, we have $d\in D(s,r)$.
Let $a'=(a'_1,\ldots,a'_{r-1}) := \delta_{r-1}^{-1}(d') \in A(s,r-1)$, and define $a = (a_1,\ldots, a_r):=\delta_{r}^{-1}(d) \in A(s,r)$.
It is straightforward to check that $a' = \partial_{u,v}(a)$. By condition (ii), we can thus apply the Tableau Reduction Lemma \ref{tablred}, which yields
\[
h_{p,r}(a) - h_{p,r-1}(a') = \frac{1}{p+1}\left(p+\frac{1}{p^{2(v-u-1)}}\right) \left(\frac{U}{p^{a_{u+1}}}+p^{a_v} V\right)
+ \frac{1}{p^2-1}\left(1- \frac{1}{p^{2(v-u-1)}}\right),
\]
where
\[
U= \sum_{k=1}^u p^{a_k} <  \sum_{k=0}^{\infty} p^{a_u-k} =\frac{p^{1+a_u}}{p-1}
\]
and  
\[
V= \sum_{i=v+1}^r \frac{1}{p^{a_i}} < \sum_{i=0}^{\infty} \frac{1}{p^{a_{v+1}+i}} = \frac{p^{1-a_{v+1}}}{p-1}\,. 
\]
Putting this together and applying $a_{u+1}-a_u =d_u =d'_u-1 =2$ and 
$a_{v+1}-a_v = d_v =d'_v-1 = 2$ by virtue of (i), we have
\begin{align*}  
\begin{split}
h_{p,r}(a) - h_{p,r-1}(a') &< \frac{1}{p^2-1}\left(p+\frac{1}{p^{2(v-u-1)}}\right) \left(\frac{p^{a_u+1}}{p^{a_{u+1}}}+\frac{p^{a_v}}{p^{a_{v+1}-1}}\right) + \frac{1}{p^2-1}\left(1- \frac{1}{p^{2(v-u-1)}}\right)\\
&= \frac{2}{p(p^2-1)}\left(p+\frac{1}{p^{2(v-u-1)}}\right) + \frac{1}{p^2-1}\left(1- \frac{1}{p^{2(v-u-1)}}\right) \\
&= \frac{1}{p^2-1}\left(3-\frac{p-2}{p^{2v-2u-1}}\right)\\
&\leq \frac{1}{p-1}
\end{split}
\end{align*}
for all primes $p$.
By Lemma \ref{hilfssatz}, this 
concludes the proof.

\qed

\begin{proposition}{\label{2221}}
Let $p$ be a prime. 
If $d=(d_1,\ldots,d_{r-1})\in D(s,r)$ for some integers $4\leq r \leq s$ with $d_u=1$ for some $u$  satisfying $2\leq u\leq r-2$ and $d_j=2$ for $j\neq u$ , then $\Ene{{\cal D}(\delta_{r}^{-1}(d))}{p^s} < \Ene{{\cal D}(\delta_{r}^{-1}(d'))}{p^s}$, where $d':=(2,2,\ldots,2,2,1)\in D(s,r)$.
\end{proposition}
 
{\sc Proof. }
This is a special case of \cite{SA3}, Proposition 3.2. Instead of becoming acquainted with the notation there, the reader might be well advised to look at the involved delta tableaux and do some calculations similar to those above.

\qed

%
%

\section{Proof of the Theorem}

As mentioned at the end of the introduction, the proof of Theorem \ref{thm} uses an accordion-like compression and expansion procedure. Starting with an arbitrary $a^{(0)}\in A(s):= \bigcup_{r=2}^s A(s,r)$, a finite sequence $(a^{(\ell)})_{0\leq \ell \leq m}$, say,  of admissible $a^{(\ell)} \in A(s)$  satisfying 
\begin{equation}  \label{seq1}
\Ene{{\cal D}(a^{(0)})}{p^s} < \Ene{{\cal D}(a^{(1)})}{p^s} < \ldots < \Ene{{\cal D}(a^{(m)})}{p^s},
\end{equation}
can be generated,
where $a^{(m)}$ is one of the maximising exponent tuples to be found in Theorem \ref{thm}.
It simplifies the matter if we perform our compression and expansion procedure on the delta vectors $d^{(\ell)}:=\delta_{r(\ell)}(a^{(\ell)}) \in D(s,r(\ell))$, $0\leq \ell \leq m$, for suitable lenghts parameters $r(\ell)$. By this (\ref{seq1}) is replaced by
\begin{equation}  \label{seq2}
\Ene{{\cal D}(\delta^{-1}_{r(0)}(d^{(0)}))}{p^s} < \Ene{{\cal D}(\delta^{-1}_{r(1)}(d^{(1)}))}{p^s} < \ldots < \Ene{{\cal D}(\delta^{-1}_{r(m)}(d^{(m)}))}{p^s}.
\end{equation}
Each transformation $d^{(\ell)} \longrightarrow d^{(\ell+1)}$ corresponds to the application of one of the Propositions in Section \ref{Compare}. In order to indicate which Proposition is used and which effect it has, we label the transformation arrows according to the following table:
\[
\begin{array}{|c|l||l|l|} \hline
\mbox{Transf.} & \mbox{Prop.}&\mbox{Effect}&\mbox{Condition} \\\hline
   & \\[-.21in]
\stackrel{\mbox{\tiny\rm Ia}}{\longrightarrow} &  \mbox{\ref{maxd=3}(i)}&
\mbox{\small Replaces entry $d_u\geq 4$ by successive entries $2, d_u-2$}& ./.\\
\stackrel{\mbox{\tiny\rm Ib}}{\longrightarrow} &   \mbox{\ref{maxd=3}(ii)}&
\mbox{\small Replaces $d_u\geq 3$ by successive entries $2, d_u-2$}&\mbox{All $d_j\geq 2$}  \\
\stackrel{\mbox{\tiny\rm II}}{\longrightarrow} &   \mbox{\ref{maxd=3mindd=1}}&
\mbox{\small Repl. $d_u=1,d_v=3$ or $d_u=3,d_v=1$ ($u<v$) by $2,2$}&\mbox{\small $d_j=2$ for $u < j <v$}  \\
\stackrel{\mbox{\tiny\rm III}}{\longrightarrow} &   \mbox{\ref{twice1} }& 
\mbox{\small Replaces $d_u=d_v=1$ ($u<v$) by $2$}&\mbox{\small $d_j=2$ for $u < j <v$} \\
\stackrel{\mbox{\tiny\rm IV}}{\longrightarrow} &   \mbox{\ref{twice3}} & 
\mbox{\small Replaces $d_u=d_v=3$ ($u<v$) by $2,2,2$}&\mbox{\small $d_j=2$ for $u < j <v$} \\
\stackrel{\mbox{\tiny\rm V}}{\longrightarrow} &   \mbox{\ref{2221}}   & 
\mbox{\small Shifts $d_u=1$ ($u<r-1$) to the right}&\mbox{\small $d_j=2$ for $j\neq uv$}\\ \hline
\end{array}
\]
The following proof of Theorem \ref{thm} demonstrates in which order the transformations can be applied to an initial delta vector $d^{(0)}$. Observe that the sequence (\ref{seq2}) is by no means unique, but the final delta vector $d^{(m)}$ very well is (with the solitary ambiguity for $p=2$ mentioned in Theorem \ref{thm}(i)). Example \ref{exseq} right behind the proof will illustrate the transformation process.

{\sc Proof of Theorem \ref{thm}. } \newline
In \cite{SA1}, Theorem 3.2, the maximal energies $\Emax{p^s}$ along with the corresponding divisor sets were determined for $s\leq 4$. Hence we may assume $s\geq 5$.
Let $a$ be an exponent tuple having the property that ${\cal D}(a)$ maximises the energy of $\Icg{\cal D}{p^s}$, i.e. $\Ene{{\cal D}(a)}{p^s} = \Emax{p^s}$. In Section \ref{intro} we saw that necessarily  $a\in A(s,r)$ for some $r$ satisfying $2\leq r \leq s$. Define $d=(d_1,\ldots,d_{r-1}):= \delta_r(a)\in D(s,r)$.

It follows from Proposition \ref{maxd=3}(i) that $d_j \in \{1,2,3\}$ for $1\leq j \leq r-1$. Thus we can visualise $d$ as an $(r-1)$-tuple possibly containing some entries $1$ and $3$ as well as entries $2$ filling the gaps. If $d$ contained two or more entries from the set $\{1,3\}$, this would contradict at least one of the Propositions \ref{maxd=3mindd=1}, \ref{twice1} or \ref{twice3}, except for the special case $d=(1,2,2,\ldots,2,1)$ for $p=2$ from Proposition \ref{twice1}. Therefore, only three situations remain: (i) $d=(2,2,\ldots,2)$,
(ii) $d=(1,2,2,\ldots,2,1)$ and $p=2$, (iii) $d$ has a single entry $1$ or $3$ with all other entries equal to $2$.  


\underline{Case 1}: $s$ is odd.
\newline
Since $\sum_{j=1}^{r-1} d_j = s-1$ is even, we have $d=(2,2,\ldots,2)\in D(s,\frac{s+1}{2})$ or, only in case $p=2$,  $d=(1,2,2,\ldots,2,1)\in D(s,\frac{s+3}{2})$. By use of (\ref{ft3}) and comparison of the delta tableaux corresponding to $\delta_{\frac{s+1}{2}}^{-1}(2,2,\ldots,2))$ and  $\delta_{\frac{s+3}{2}}^{-1}((1,2,2,\ldots,2,1))$, it is an easy exercise to check that 
\[
\Ene{{\cal D}(\delta_{\frac{s+1}{2}}^{-1}(2,2,\ldots,2))}{2^s} = \Ene{{\cal D}(\delta_{\frac{s+3}{2}}^{-1}(1,2,2,\ldots,2,1))}{2^s}.
\]
It remains to prove the formula for $\Emax{p^s}$. We have 
\begin{align}   \label{hpodd}
\begin{split}
h_{p,\frac{s+1}{2}}(0,2,4,\ldots, s-3,s-1) &= \sum_{k=1}^{\frac{s-1}{2}} \sum_{i=k+1}^{\frac{s+1}{2}} \frac{1}{p^{2(i-k)}} = 
\sum_{j=1}^{\frac{s-1}{2}} \frac{1}{p^{2j}}\left(\frac{s+1}{2}-j\right)   \\
&= \frac{(s-1)p^{s+1}-(s+1)p^{s-1}+2}{2(p^2-1)^2 p^{s-1}}\, ,
\end{split}
\end{align} 
hence by (\ref{ft3})
\begin{align*}
\Emax{p^s} &= \Ene{{\cal D}((0,2,4,\ldots,s-3,s-1))}{p^s} \\
&= 2(p-1)p^{s-1}\left( \frac{s+1}{2}- (p-1)h_{p,\frac{s+1}{2}}(0,2,4,\ldots, s-3,s-1)\right)\\
&= \frac{1}{(p+1)^2}\left((s+1)(p^2-1)p^s + 2(p^{s+1}-1)\right),
\end{align*}
which proves (i).

\underline{Case 2}: $s$ is even.
\newline
Since $\sum_{j=1}^{r-1} d_j = s-1$ is odd, $d$ must have a single entry $1$ or $3$ with all other entries equal to $2$. By Proposition \ref{maxd=3}(ii) the unique odd entry is $1$, and by
Proposition \ref{2221} and the symmetry property (\ref{symmhp}), we necessarily have
$d=(2,2,\ldots,2,2,1)\in D(s,\frac{s+2}{2}))$ or $d=(1,2,2,\ldots,2,2)\in D(s,\frac{s+2}{2}))$ with equal corresponding energy.
With regard to $\Emax{p^s}$, we have in this case
\begin{align}  \label{hpeven}
\begin{split}
h_{p,\frac{s+2}{2}}(0,2,4,\ldots, s-2,s-1) &= \sum_{k=1}^{\frac{s-2}{2}} \sum_{i=k+1}^{\frac{s}{2}} \frac{1}{p^{2(i-k)}} + \sum_{k=1}^{\frac{s}{2}} \frac{1}{p^{s-1-2(k-1)}}\\
&= 
\sum_{j=1}^{\frac{s-2}{2}} \frac{1}{p^{2j}}\left(\frac{s}{2}-j\right)  + \frac{1}{p^{s-1}}\sum_{k=0}^{\frac{s-2}{2}} p^{2k} \\
&=   \frac{(s-2)p^s-sp^{s-2}+2}{2(p^2-1)^2 p^{s-2}}+  \frac{p^s-1}{(p^2-1)p^{s-1}}\, ,
\end{split}
\end{align} 
hence by (\ref{ft3})
\begin{align*}
\Emax{p^s} &= \Ene{{\cal D}((0,2,4,\ldots,s-2,s-1))}{p^s} \\
&= 2(p-1)p^{s-1}\left( \frac{s+2}{2}- (p-1)h_{p,\frac{s+2}{2}}(0,2,4,\ldots, s-2,s-1)\right)\\
&= \frac{1}{(p+1)^2}\left(s(p^2-1)p^s + 2(2p^{s+1}-p^{s-1}+ p^2-p-1)\right),
\end{align*}
and this proves (ii).


\qed

%
%
%
%
 
\begin{example}   \label{exseq}
\end{example} 
Let $a^{(0)}:= (0,5,6,9,12,14,15,16,22,23,24,27,29) \in A(30,13)$, hence
\[
d^{(0)}:= \delta_{13}(a^{(0)})=(5,1,3,3,2,1,1,6,1,1,3,2) \in D(30,13).
\] 
Step by step, we obtain
{\small
\[
\begin{array}{|c||ll|c|r|r|}
\hline
\ell & \multicolumn{1}{c}{d^{(\ell)}} &&r(\ell) & \Ene{{\cal D}(\delta_{r(\ell)}^{-1}(d^{(\ell)}))}{2^{30}}&
  \Ene{{\cal D}(\delta_{r(\ell)}^{-1}(d^{(\ell)}))}{3^{30}}\\
\hline
0 & (\underline{5},1,3,3,2,1,1,6,1,1,3,2) & & 13 &  9~167~691~382 & 2~293~430~091~118~444  \\[-.1in]
& &\big\downarrow{^{\mbox{\!\tiny{Ia}}}}&& & \\[-.08in]
1 & (\mbox{\bf 2},\mbox{\bf 3},1,3,3,2,1,1,\underline{6},1,1,3,2) && 14 &  9~761~773~390 &2~479~571~746~112~800\\[-.1in]
& &\big\downarrow{^{\mbox{\!\tiny{Ia}}}}&&  &\\[-.08in]
2 & (2,3,1,3,3,2,1,1,\mbox{\bf 2},\mbox{\bf \underline{4}},1,1,3,2) && 15 &  10~226~403~150& 2~655~370~924~580~476\\[-.1in]
& &\big\downarrow{^{\mbox{\!\tiny{Ia}}}}&&  &\\[-.08in]
3 & (2,3,1,3,3,2,1,\underline{1},2,\mbox{\bf 2},\mbox{\bf 2},\underline{1},1,3,2) && 16 & 10~429~199~182  &2~770~612~868~608~768  \\[-.1in]
& &\big\downarrow{^{\mbox{\!\tiny{III}}}}&& & \\[-.08in]
4 & (2,3,\underline{1},\underline{3},3,2,1,\mbox{\bf 2},2,2,2,1,3,2) & & 15 &  10~869~926~478 &2~937~991~189~453~948 \\[-.1in]
& &\big\downarrow{^{\mbox{\!\tiny{II}}}}&&  &\\[-.08in]
5 & (2,3,\mbox{\bf 2},\mbox{\bf 2},3,2,\underline{1},2,2,2,2,\underline{1},3,2) & & 15 &  11~022~317~518 &3~022~615~444~978~108 \\[-.1in]
& &\big\downarrow{^{\mbox{\!\tiny{III}}}}&& & \\[-.08in]
6 & (2,3,2,2,\underline{3},2,\mbox{\bf 2},2,2,2,2,\underline{3},2) & & 14 &  11~182~822~222 &3~112~785~070~640~560\\[-.1in]
& &\big\downarrow{^{\mbox{\!\tiny{IV}}}}&&  &\\[-.08in]
7 & (2,\underline{3},2,2,\mbox{\bf 2},2,2,2,2,2,2,\mbox{\bf 2},\mbox{\bf 2},2) & & 15 & 11~438~333~038 &3~216~413~472~521~788 \\[-.1in]
& &\big\downarrow{^{\mbox{\!\tiny{V}}}}&& & \\[-.08in]
8 & (2,\mbox{\bf \underline{1}},\mbox{\bf 2},2,2,2,2,2,2,2,2,2,2,2,2) & & 16 & 11~483~072~286  &3~218~955~338~350~144\\[-.1in]
& &\big\downarrow{^{\mbox{\!\tiny{Ib}}}}&& & \\[-.08in]
9 & (2,2,2,2,2,2,2,2,2,2,2,2,2,2,\mbox{\bf 1}) & & 16 &  11~572~550~770 &3~234~206~533~320~112\\ \hline
\end{array}
\]}
According to Theorem \ref{thm}(ii)~ we have
\[ \Emax{2^{30}} = \Ene{{\cal D}((0,2,4,\ldots,28,29))}{2^{30}} = 11~572~550~770
\]
and
\[ \Emax{3^{30}} = \Ene{{\cal D}((0,2,4,\ldots,28,29))}{3^{30}} = 3~234~206~533~320~112.
\]

%
%


\begin{thebibliography}{999}
\addcontentsline{toc}{}{References}


\bibitem{AHM}
O. Ahmadi and N. Alon and I.F. Blake and I.E. Shparlinski,
Graphs with integral spectrum,
Linear Algebra Appl. {\bf 430} (2009), 547-552.

\bibitem{AKH} R. Akhtar and M. Boggess and T. Jackson-Henderson and I. Jim{\'e}nez and R. Karpman and A. Kinzel and D. Pritikin,
On the unitary Cayley graph of a finite ring,   
Electron. J. Combin. {\bf 16} (2009), Research Paper R117, 13 pp. (electronic). 

\bibitem{BAP} R.B. Bapat and S. Pati,
Energy of a graph is never an odd integer,
Bull. Kerala Math. Assoc. {\bf 1} (2004), 129-132.

\bibitem{BAS} M. Ba{\v s}i{\'c} and A. Ili{\'c},
{On the clique number of integral circulant graphs},
Appl. Math. Lett. {\bf 22} (2009), 1406-1411.

\bibitem{BAS2} M. Ba{\v s}i{\'c} and A. Ili{\'c},
{On the Automorphism Group of Integral Circulant Graphs},
Electron. J. Combin. {\bf 18} (2011), Research Paper P68, 13 pp. (electronic). 

\bibitem{BAS3} M. Ba\v{s}i\'c and M.D. Petkovi\'c,
{Perfect state transfer in integral circulant graphs of non-square-free order},
Linear Algebra Appl. {\bf 433} (2010), 149-163

\bibitem{BAS4} M. Ba\v si\'c and M.D. Petkovi\'c and D. Stevanovi\'c,
{Perfect state transfer in integral circulant graphs},
Appl. Math. Lett. {\bf 22} (2009), 1117-1121

\bibitem{BEA} N. de Beaudrap,
{On restricted unitary {C}ayley graphs and symplectic transformations modulo {$n$}},
Electron. J. Combin. {\bf 17} (2010), Research Paper R69, 26 pp. (electronic). 

\bibitem{BER} P. Berrizbeitia and R. E. Giudici,
{On cycles in the sequence of unitary Cayley graphs},
Discrete Math. {\bf 282} (2004), 239-243.

\bibitem{BRU} R.A. Brualdi,
Energy of a graph, AIM Workshop Notes, 2006.

\bibitem{DAV} P. J. Davis,
{Circulant matrices},
John Wiley \&\ Sons, New York-Chichester-Brisbane, 1979.
  
\bibitem{DEJ} I. J. Dejter and R. E. Giudici,
{On unitary Cayley graphs},
J. Combin. Math. Combin. Comput. {\bf 18} (1995), 121-124.
   
\bibitem{DRO} A. Droll,
{A classification of {R}amanujan unitary {C}ayley graphs},
Electron. J. Combin. {\bf 17} (2010), Research Note N29, 6 pp. (electronic). 

\bibitem{GUT} I. Gutman,
{The energy of a graph},
Ber. Math.-Stat. Sekt. Forschungszent. Graz {\bf 103}, 1978.
   
\bibitem{ILI} A. Ili\'{c},
The energy of unitary Cayley graphs,
Linear Algebra Appl. {\bf 431} (2009), 1881-1889.

\bibitem{KIA} D. Kiani and M.M.H. Aghaei and Y. Meemark and B. Suntornpoch,
The energy of unitary Cayley graphs and gcd-graphs,
Linear Algebra Appl. {\bf 435} (2011), 1336-1343.

\bibitem{KLO} W. Klotz and T. Sander,
Some properties of unitary Cayley graphs,
Electron. J. Combin. {\bf 14} (2007), Research Paper R45, 12 pp. (electronic).

\bibitem{KLO2} W. Klotz and T. Sander,
{Integral {C}ayley graphs over abelian groups},
Electron. J. Combin. {\bf 17} (2010), Research Paper R81, 13 pp. (electronic). 

\bibitem{KOO} J.H. Koolen and V. Moulton,
{Maximal energy graphs},
Adv. Appl. Math. {\bf 26} (2001), 47-52.

\bibitem{PET} M.D. Petkovi\'c and M. Ba\v{s}i\'c,
{Further results on the perfect state transfer in integral circulant graphs},
Comput. Math. Appl. {\bf 61} (2011), 300-312

\bibitem{PIR} S. Pirzada and I. Gutman,
Energy of a graph is never the square root of an odd integer,
Appl. Analysis and Discr. Math. {\bf 2} (2008), 118-121.

\bibitem{RAM} H.N. Ramaswamy and C.R. Veena,
{On the Energy of Unitary Cayley Graphs},
Electron. J. Combin. {\bf 16} (2009), Research Note N24, 8 pp. (electronic).

\bibitem{SA1} J.W. Sander and T. Sander,
{The energy of integral circulant graphs with prime power order},
Appl. Anal. Discrete Math. {\bf 5} (2011), 22-36.

\bibitem{SA2} J.W. Sander and T. Sander,
{Integral circulant graphs of prime power order with maximal energy},
Linear Algebra Appl. {\bf 435} (2011), 3212-3232.    

\bibitem{SA3} J.W. Sander and T. Sander,
{The maximal energy of classes of integral circulant graphs},
submitted.

\bibitem{SAX} N. Saxena and S. Severini and I.E. Shparlinski,
{Parameters of integral circulant graphs and periodic quantum dynamics},
Int. J. Quantum Inf. {\bf 5} (2007), 417-430.        

\bibitem{SHP} I. Shparlinski,
{On the energy of some circulant graphs},
Linear Algebra Appl. {\bf 414} (2006), 378-382.

\bibitem{SO} W. So,
Integral circulant graphs,
Discrete Math. {\bf 306} (2005), 153-158.

\end{thebibliography}
\end{document}